\documentclass[a4paper,12pt]{amsart}
\usepackage{amsmath,amssymb,amsfonts,graphicx,color,fancyhdr}
\usepackage[latin1]{inputenc}
\usepackage{psfrag}

\newtheorem{theorem}{Theorem}[section]
\newtheorem{proposition}[theorem]{Proposition}
\newtheorem{lemma}[theorem]{Lemma}
\newtheorem{definition}[theorem]{Definition}

\newtheorem{remark}[theorem]{Remark}
\newtheorem{corollary}[theorem]{Corollary}
\newtheorem{example}[theorem]{Example}

\usepackage[colorlinks=true, pdfstartview=FitV, linkcolor=blue, citecolor=blue, urlcolor=blue]{hyperref}
\usepackage{color}
\usepackage{psfrag,caption}

\def\be#1 {\begin{equation} \label{#1}}
\newcommand{\ee}{\end{equation}}

\def\sqw{\hbox{\rlap{\leavevmode\raise.3ex\hbox{$\sqcap$}}$%
\sqcup$}}
\def\findem{\ifmmode\sqw\else{\ifhmode\unskip\fi\nobreak\hfil
\penalty50\hskip1em\null\nobreak\hfil\sqw
\parfillskip=0pt\finalhyphendemerits=0\endgraf}\fi}

\newcommand{\R}{{\mathbb {R}}}

\newcommand{\Z}{{\mathbb Z}}

\begin{document}

\title
[Fourier Multipliers and Littlewood-Paley For Modulation Spaces]
{\bf Fourier Multipliers and Littlewood-Paley For Modulation Spaces}

\author{Parasar Mohanty}
\address{Department of Mathematics \\
Indian Institute of Technology Kanpur-208016, India}
\email{parasar@iitk.ac.in}

\author{Saurabh Shrivastava}
\address{Department of Mathematics \\
Indian Institute of Science Education and Research Bhopal, India}

\email{saurabhk@iiserb.ac.in}
\begin{abstract}
In this paper we have studied Fourier multipliers and Littlewood-Paley square functions in the context of modulation 
spaces. We have also proved that any bounded linear operator from modulation space $\mathcal{M}_{p,q}(\R^n),~1\leq p,q\leq \infty,~$ 
into itself possesses an $l_2-$valued extension. This is an analogue of a well known result due to Marcinkiewicz and Zygmund on classical 
$L^p-$spaces.
\end{abstract}

\date{\today}

\maketitle
\section{\bf Introduction}

The theory of modulation spaces has been developed substantially in the last decade. Modulation spaces provide quantitative information
about time-frequency concentration of functions and distributions. The modulation spaces are defined in terms of short-time Fourier transform.
The short-time Fourier transform of a function is defined as inner product of the function with respect to a time-frequency shift of another function,
known as a window function (for precise definition see Section~\ref{pre}). Modulation spaces have found their usefulness in applications as well
as in pure mathematics. They play an useful role in the theory of pseudo-differential operators. We refer the interested reader to~\cite{T1,T2}
for this connection.

The purpose of this paper is to study Fourier multipliers and Littlewood-Paley operators in the context of modulation spaces.

The paper is organized as follows~:~

In Section~\ref{pre} we set notation and give a brief introduction to modulation spaces.
We will discuss various equivalent definitions of modulation spaces. We will also mention some
basic properties of short-time Fourier transform and modulation spaces in this section.

In Section~\ref{fm} we will address some  natural
questions about Fourier multipliers on modulation spaces. It is a well known fact that classical
$L^p-$multipliers are always Fourier multipliers for
respective modulation spaces $\mathcal{M}_{p,q}(\R^n),~1\leq p,q\leq \infty.$ So, the real interest here is to
investigate if non $L^p-$multipliers give rise to Fourier multipliers for modulation spaces.
We will discuss some known results~\cite{BGGO,BGOR,FN} in this context and
also provide some new examples of non $L^p-$multipliers, which become multipliers on modulation spaces.
In this work we have generalized the main result of~\cite{BGGO}. The next question is to find bounded measurable functions, which are not Fourier multipliers for
modulation spaces $\mathcal{M}_{p,q}(\R^n),~p\neq 2.$ We will show that some known classical $L^p$ results
provide examples of such functions for modulation spaces as well. Moreover, we will construct a new example in this
direction, which will be relevant to answer some questions about Littlewood-Paley operators as well. Next, we will study inclusion relations between Fourier
multiplier spaces on modulation spaces.

Section~\ref{lp} is devoted to the theory of Littlewood-Paley square functions on modulation spaces.

In Section~\ref{mz} we establish the modulation space analogue of a classical result due to Marcinkiewicz and Zygmund
about $l_2-$valued extension of bounded linear operators.

\section{\bf Notation and Preliminaries}\label{pre}

Let $f$ be a complex-valued function defined on $\R^n.$ Consider the operations of translation, modulation,
and dilation defined as follows:
\begin{itemize}
 \item Translation operator~:~$\tau_yf(x)= f(x-y),~ {\rm for}\;x,y\in \R^n.$
\item Modulation operator~:~$M_{\xi}f(x)= e^{2 \pi i x.\xi}f(x),~ {\rm for}\;x,\xi\in \R^n.$
\item Dilation~:~$D_{\lambda}f(x)= f(\lambda x),~ {\rm for}\;\lambda>0;~x\in \R^n.$
\end{itemize}
We denote $\mathcal{S}(\R^n)$ the space of Schwartz class functions on $\R^n$ and $\mathcal{S'}(\R^n)$ the space of
tempered distributions on $\R^n.$
\medskip

We use the following definition for the Fourier transform of an $L^1(\R^n)$ function:
$$\hat{f}(\xi)= \int_{\R^n} f(x) e^{-2 \pi i x.\xi} dx,~\xi\in \R^n.$$

\begin{definition}[Short-time Fourier Transform] Let $g\in \mathcal{S}(\R^n)$ be a non-zero function.
The short-time Fourier transform of $f\in \mathcal{S'}(\R^n)$ with respect to $g$ is defined as:
\begin{eqnarray}
V_g(f)(x,\xi)&=& \left<f, M_{\xi} \tau_x g \right>.
\end{eqnarray}
Observe that if in addition $f$ is a nice function, then
\begin{eqnarray}
V_g(f)(x,\xi)&=& \int_{\R^n} e^{-2 \pi i t.\xi}~ \overline{g(t-x)} f(t) dt.
\end{eqnarray}
We also say that $V_g(f)$ is the short-time Fourier transform of $f$ with respect to the window $g.$
\end{definition}

For $x,\xi\in \R^n,~M_{\xi} \tau_x g$ is said to be the time-frequency shift of $g$ by $(x,\xi).$  Thus
the short-time Fourier transform $V_g f$ is the inner product of $f$ with respect to time-frequency shift
of $g.$ It is also interpreted as a simultaneous time-frequency representation of $f$ in the time-frequency plane.
Moreover, it also occurs under the names of ``cross ambiguity function'' and ``cross Wigner distribution'' with slightly
different formulations. For more details about cross ambiguity function and cross Wigner distribution
see~Chapter~$4$~in~\cite{G}. A few  different forms of $V_gf$ are listed as below.
\begin{lemma}\label{stft1}\cite{G} If $f,g\in L^2(\R^n),$ then $V_gf$ is uniformly continuous function on $\R^{2n},$ and
\begin{eqnarray*}
V_g(f)(x,\xi)&=& \left<f, M_{\xi} \tau_x g \right>\\
&=& \left<\hat f, \tau_{\xi} M_{-x}\hat g \right>\\
&=& \widehat{(f.\tau_x \bar{g})}(\xi)\\
&=& e^{-2 \pi i x.\xi}\widehat{(\hat f.\tau_{\xi} \bar{\hat g})}(-x)\\
&=& e^{-2 \pi i x.\xi} V_{\hat g}(\hat f)(\xi,-x)\\
&=& e^{-2 \pi i x.\xi} (f\ast M_{\xi} g^*)(x)\\
&=& (\hat f\ast M_{-x} \hat g^*)(\xi)\\
&=& e^{- \pi i x.\xi}\int_{\R^n} f(t+\frac{x}{2})\overline{g(t-\frac{x}{2})}e^{-2 \pi i t.\xi}dt,
\end{eqnarray*}
where $g^*(x)=\overline {g(-x)}.$
\end{lemma}
The proof of the above lemma follows easily using Parseval's formula. For $f,g\in L^2(\R^n),$ the quantity
$\int_{\R^n} f(t+\frac{x}{2})\bar g(t-\frac{x}{2})e^{-2 \pi i t.\xi}dt$ is referred to as cross ambiguity
function of $f$ and $g.$
\medskip

The short-time Fourier transform enjoys some similar properties like the classical Fourier transform. Some of
them are listed in the following proposition.
\begin{proposition}\label{stft2}\cite{G} We have
\begin{enumerate}
\item {\bf Orthogonality property~:}~$\left<V_{g_1}(f_1), V_{g_2}(f_2)\right> = \left<f_1, f_2\right>
\overline{\left<g_1, g_2\right>},~f_1, f_2, g_1, g_2 \in L^2(\R^n).$
\item {\bf $L^2-$norm relation~:}~$\|V_gf\|_{L^2(\R^{2n})}= \|f\|_{L^2(\R^n)} \|g\|_{L^2(\R^n)},~f, g \in L^2(\R^n).$
\item {\bf Inversion formula~:}~If $g_1, g_2 \in L^2(\R^n)$ are such that $\left<g_1, g_2\right>\neq 0,$
then for all $f\in L^2(\R^n)$ we have
$$f=\frac{1}{\left<g_1, g_2\right>}\int_{R^{2n}}V_{g_1}(f)(x,\xi)M_{\xi} \tau_x g_2 d\xi dx.$$
\end{enumerate}
\end{proposition}
Properties $(2)$ and $(3)$ follow immediately once we have property $(1)$ and property $(1)$ can be proved
by using Parseval's formula.
\medskip

We shall now define modulation spaces.
\begin{definition}[Modulation Spaces]\cite{G} Let $1\leq p,q \leq \infty$ and let $g$ be a window function. Then the
modulation space $\mathcal{M}_{p,q}(\R^n)$ is the space of all tempered distributions $f\in \mathcal{S'}(\R^n)$
for which the following mixed norm is finite:
\begin{eqnarray}
\|f\|_{\mathcal{M}_{p,q}}= \left(\int_{\R^n} \left(\int_{\R^n} |V_g(f)(x,\xi)|^p dx\right)^{\frac{q}{p}}d\xi\right)^{\frac{1}{q}},
\end{eqnarray}
with the usual modifications when $p$ and/or $q$ are infinite.
\end{definition}
The definition of modulation spaces is independent of choice of the window function $g$ in the sense of equivalent
norms. We use the notation $\mathcal{M}_p(\R^n)$ for modulation space $\mathcal{M}_{p,p}(\R^n).$
\medskip

The modulation spaces have various equivalent definitions. For example, the norm $\|.\|_{\mathcal{M}_{p,q}}$
has the following equivalent formulation (also see~\cite{FN}).

Let $\phi$ be a smooth function defined on $\R^n$ such that ${\it supp}~ \phi \subseteq [-1,1]^n$ and
$\sum_{k\in \Z^n}\phi(\xi-k)=1$ for all $\xi\in \R^n.$ Set $\phi_k(\xi)= \phi(\xi-k)$ and let $T_{\phi_k}$
be the Fourier multiplier operator given by $\widehat{T_{\phi_k}f}(\xi)= \phi_k(\xi) \hat{f}(\xi).$ Then, we have
\begin{eqnarray}\label{def2}
\|f\|_{\mathcal{M}_{p,q}} \simeq \left(\sum_{k\in \Z^n}\|T_{\phi_k}f\|_p^q\right)^{\frac{1}{q}},
\end{eqnarray}
with appropriate modification when $p$ and/or $q$ are infinite. Here the notation $A\simeq B$ means that there
are two positive constants $C_1,C_2$ such that $C_1 A\leq B \leq C_2 A.$

The above definition turns out to be very useful in order to study Fourier multipliers on modulation spaces.

Next, we present yet another definition of modulation spaces. This is given via Gabor frames and plays
a key role in order to study simultaneously the local time and frequency behaviour of functions. Let us first define Gabor frames.
\begin{definition}[Gabor frames]\cite{G,Oko}
Let $\phi\in L^2(\R^n)$ and $\alpha, \beta>0$. Then we say that $\{M_{\alpha l}\tau_{\beta k}\phi\}_{k,l\in \Z^n}$
is a Gabor frame for $L^2(\R^n)$ if
\begin{eqnarray}
\|f\|^2_{L^2(\R^n)}\simeq \sum_{k,l\in \Z^n}|\left<f, M_{\alpha l}\tau_{\beta k}\phi \right>|^2.
\end{eqnarray}
\end{definition}
The spaces $\mathcal{M}_{p,q}(\R^n)$ are characterized by means of Gabor frames in the form of the following theorem:
\begin{theorem} \cite{BGGO,G} Let $\phi\in \mathcal{M}_1(\R^n)$ be such that $\{M_{\alpha l}\tau_{\beta k}\phi\}_{k,l\in \Z^n}$
is a Gabor frame for $L^2(\R^n).$ Then for $1\leq p,q\leq \infty,$ there exists $\psi\in \mathcal{M}_1(\R^n)$ (dual frame)
such that every $f\in \mathcal{M}_{p,q}(\R^n)$ has a Gabor expansion that converges unconditionally (or $\it{weak}^*$
unconditionally when $p=\infty$ or $q=\infty$), i.e., we have
\begin{eqnarray}
f = \sum_{k,l\in \Z^n}\left<f, M_{\alpha l}\tau_{\beta k}\psi \right> M_{\alpha l}\tau_{\beta k}\phi,~ \forall f\in \mathcal{M}_{p,q}(\R^n).
\end{eqnarray}
Moreover, the following equivalence holds:
\begin{eqnarray*}
\|f\|_{\mathcal{M}_{p,q}}
&\simeq& \left(\sum_{l\in \Z^n}\left(\sum_{k\in \Z^n}|\left<f, M_{\alpha l}\tau_{\beta k}\phi\right>|^p\right)^{\frac{q}{p}}\right)^{\frac{1}{q}}\\
&\simeq& \left(\sum_{l\in \Z^n}\left(\sum_{k\in \Z^n}|\left<f, M_{\alpha l}\tau_{\beta k}\psi\right>|^p\right)^{\frac{q}{p}}\right)^{\frac{1}{q}}.
\end{eqnarray*}
\end{theorem}
We shall use any of these definitions for modulation spaces as per our requirement. We list here some basic and important
facts about modulation spaces, whose proofs are available in various literature~\cite{G,ST}.
\begin{proposition}\label{prop} The following are some of the important properties of modulation spaces~:
\begin{enumerate}
\item The space of Schwartz class functions $\mathcal{S}(\R^n)$ is dense in $\mathcal{M}_{p,q}(\R^n)$ for all
$1\leq p,q<\infty.$
\item The modulation spaces are invariant under the operations of translation, modulation, and dilation.
\item The dual of $\mathcal{M}_{p,q}(\R^n),~1\leq p,q<\infty,$ is $\mathcal{M}_{p',q'}(\R^n),$ where
$\frac{1}{p}+\frac{1}{p'}=1$ and $\frac{1}{q}+\frac{1}{q'}=1.$
\item The modulation space $\mathcal{M}_1(\R^n)$ becomes a Banach algebra under both pointwise multiplication
and convolution. Moreover, it is invariant under the Fourier transform. Also, $\mathcal{M}_1(\R^n)$ is called
the Feichtinger algebra.
\item $\mathcal{M}_2(\R^n)=L^2(\R^n).$
\item $ \mathcal{M}_{p_1,q_1}(\R^n)\subseteq \mathcal{M}_{p_2,q_2}(\R^n)$ whenever $p_1\leq p_2$ and $q_1\leq q_2.$
\item $\mathcal{M}_p(\R^n)\subseteq L^p(\R^n)\subseteq \mathcal{M}_{p,p'}(\R^n)$ if $1\leq p\leq 2$ and
$ \mathcal{M}_{p,p'}(\R^n)\subseteq L^p(\R^n)\subseteq \mathcal{M}_p(\R^n)$ if $2\leq p\leq\infty.$
\end{enumerate}
\end{proposition}

The next theorem is due to Feichtinger~\cite{F} about complex interpolation for modulation spaces.
\begin{theorem}\cite{F,Oko} Let $1\leq p_1,q_1<\infty$ and $1\leq p_2,q_2\leq \infty,$ and $\theta\in (0,1).$
Set $\frac{1}{p}=\frac{\theta}{p_1}+\frac{1-\theta}{p_2}$ and $\frac{1}{q}=\frac{\theta}{q_1}+\frac{1-\theta}{q_2},$ then we have
$$(\mathcal{M}_{p_1,q_1}(\R^n), \mathcal{M}_{p_2,q_2}(\R^n))_{[\theta]}= \mathcal{M}_{p,q}(\R^n).$$
\end{theorem}
\section{\bf Fourier multipliers}\label{fm}
\begin{definition}[Fourier multipliers on modulation spaces] Let $m$ be a bounded measurable function defined
on $\R^n$. We say that $m$ is a Fourier multiplier on space $\mathcal{M}_{p,q}(\R^n)$ if the linear operator
$T_m,$ defined as
$$\widehat{T_mf}= m\hat{f},~f\in \mathcal{S}(\R^n),$$
extends to a bounded linear operator from $\mathcal{M}_{p,q}(\R^n)$ into itself, i.e.,
there exists a constant $C>0$ such that for all $f\in \mathcal{M}_{p,q}(\R^n),$ we have
\begin{eqnarray*}
\|T_mf\|_{\mathcal{M}_{p,q}}\leq C \|f\|_{\mathcal{M}_{p,q}}.
\end{eqnarray*}
\end{definition}
Let $M(\mathcal{M}_{p,q}(\R^n))$ denote the space of all Fourier multiplier on modulation space
$\mathcal{M}_{p,q}(\R^n).$ We will use the notation $M(L_p(\R^n))$ for the space of classical Fourier multiplier on $L^p(\R^n).$

First, we would like to describe some known results (relevant to our setting) about Fourier multipliers
in the context of modulation spaces.

Feichtinger and Narimani~\cite{FN} gave a characterization of Fourier multipliers on modulation spaces in
terms of Wiener amalgam spaces, which are defined as follows~:
\begin{definition}[Wiener amalgam spaces]\cite{FN} Let $\phi\in \mathcal{S}(\R^n)$ be such that 
${\it supp}~\phi \in [-1,1]^n$
and $\sum\limits_{k\in \Z^n}\phi_k(\xi)=1$ for all $\xi\in \R^n,$ where $\phi_k(\xi)= \phi(\xi-k).$ The Wiener amalgam space
$W(M(L_p(\R^n)), l_{\infty})$ with respect to the partition of unity $\phi_k$ is defined as the space of all
$\sigma\in \mathcal{S'}(\R^n)$ such that $\phi_k \sigma \in M(L_p(\R^n))~\forall k\in \Z^n$ with
$ {\it sup}_k \|\phi_k \sigma\|_{M(L_p)} < \infty.$

The quantity ${\it sup}_k \|\phi_k \sigma\|_{M(L_p(\R^n))}$ is the norm on $W(M(L_p(\R^n)), l_{\infty})$ and
is denoted as $\|\sigma\|_{W(M(L_p(\R^n)), l_{\infty})}.$
\end{definition}
For a more general definition of these spaces, we refer to~\cite{FN}. The next theorem due to Feichtinger and
Narimani~\cite{FN} characterizes Fourier multipliers on modulation spaces $\mathcal{M}_{p,q}(\R^n).$
\begin{theorem}~[Theorem $17(\rm i)$\cite{FN}]\label{multiplier} Let $1\leq p,q\leq \infty.$ Then we have,
$$M(\mathcal{M}_{p,q}(\R^n))=W(M(L_p(\R^n)), l_{\infty})$$
In particular, $M(\mathcal{M}_{p,q}(\R^n))$ is independent of the parameter $q.$
\end{theorem}
Fourier multipliers on modulation spaces share some properties with classical $L^p-$multipliers. We describe here some
of them as follows:
\begin{proposition}\label{prop}We have
\begin{enumerate}
\item If $\sigma \in M(\mathcal{M}_{p,q}(\R^n)),$ then for all $x,\xi\in \R^n$
and $\lambda >0,~\tau_x \sigma, M_{\xi}\sigma, D_{\lambda}\sigma \in M(\mathcal{M}_{p,q}(\R^n)).$ Moreover,
the norms of $\tau_x \sigma$ and $M_{\xi}\sigma$ are independent of $x,\xi\in \R^n.$
\item If $\sigma_1, \sigma_2 \in M(\mathcal{M}_{p,q}(\R^n))$ then $\sigma_1. \sigma_2 \in M(\mathcal{M}_{p,q}(\R^n)).$
\item $M(\mathcal{M}_{2,q}(\R^n))=M(L_2(\R^n))=L^{\infty}(\R^n).$
\item $M(L_p(\R^n))\subseteq M(\mathcal{M}_{p,q}(\R^n)).$
\end{enumerate}
\end{proposition}
When $p\neq 2,$ inclusion in property $(4)$ of the above proposition is strict. In fact, in this section we
will discuss some examples of non $L^p-$multipliers for $p\neq 2,$ which become multipliers for the corresponding
modulation space $\mathcal{M}_{p,q}(\R^n).$ This is one of the most interesting feature of Fourier multipliers on
modulation spaces. In this direction we first recall some interesting results from~\cite{BGGO,BGOR,FN}, which
are relevant to our setting.

It is well known that the function $e^{i|\xi|^2}$ is not a Fourier multiplier on $L^p(\R^n)$ unless $p=2$ (see~\cite{H}).
In contrast to this B\'{e}nyi, Gr\"{o}chenig, Okoudjou, and  Rogers~\cite{BGOR}
proved that $e^{i|\xi|^2}$ gives rise to Fourier multiplier on all modulation spaces. More precisely, they proved that
\begin{theorem}\cite{BGOR}
For $0\leq \alpha \leq 2,$ the function $e^{i|\xi|^{\alpha}}$ is a Fourier multiplier on modulation spaces
$\mathcal{M}_{p,q}(\R^n)$ for all $1\leq p,q\leq \infty$ and $n\geq 1.$
\end{theorem}

Let $\omega$ be an interval in $\R$ and let $S_{\omega}$ denote the Fourier multiplier operator given by
$\widehat {S_{\omega}f}(\xi)= \chi_{\omega}(\xi)\hat f(\xi),~f\in\mathcal{S}(\R),$
where $\chi_{\omega}$ is the characteristic function of interval $\omega.$
It is a classical fact that $\chi_{\omega}\in M(L^p(\R)),~1<p<\infty,$ which in turn implies that
$\chi_{\omega}\in M(\mathcal{M}_{p,q}(\R))$ for all $1<p<\infty,~1\leq q\leq \infty.$ Moreover, operator norm is
independent of the interval $\omega.$ However, from the classical Littlewood-Paley theory for $L^p-$spaces, we know
that function of the form  $\sum_n a_n \chi_{\omega_n}(\xi),$ where $\omega_n$ are disjoint intervals in $\R$ and $\{a_n\}$
is a bounded sequence of real numbers, may not be a multiplier on $L^p(\R),$ unless $p=2.$ In this section, we
investigate if such functions give rise to Fourier multipliers on modulation spaces.

Given a collection $\Omega =\{\omega_n:n\in \Z\}$ of disjoint intervals in $\R$ and a bounded sequence
${\bf a}=\{a_n\},$ consider $\sigma_{\Omega,{\bf a}}(\xi)= \sum_n a_n \chi_{\omega_n}(\xi).$ Let
$H_{\Omega, {\bf a}}$ denote the linear operator given by
$\widehat{H_{\Omega, {\bf a}}f}= \sigma_{\Omega,{\bf a}} \hat f.$ Our aim is to investigate that for which sequence of intervals
$\sigma_{\Omega,{\bf a}}$ gives rise to Fourier multiplier
on $\mathcal{M}_{p,q}(\R)$ for all bounded sequences ${\bf a}=\{a_n\}.$ It is again a well known classical result
that if $\Omega=\{[2^n,2^{n+1}):n\in \Z\},$ then $\sigma_{\Omega,{\bf a}}\in M(L^p(\R))$ for all $1<p<\infty.$
Hence, in the case of dyadic intervals we can easily deduce that
$\sigma_{\Omega,{\bf a}}\in M(\mathcal{M}_{p,q}(\R))$ for all $1<p<\infty,~1\leq q\leq \infty.$

Another important collection of intervals is $\Omega=\{[n,n+1)~:n\in \Z\},$ but for this collection, there exists a bounded sequence
$\{a_n\}$ such that the associated function $\sigma_{\Omega,{\bf a}}$ is not an $L^p-$multiplier unless $p=2.$
Unlike the $L^p-$case B\'{e}nyi, Grafakos, Gr\"{o}chenig, and
Okoudjou~\cite{BGGO} proved that even for this sequence $\sigma_{\Omega,{\bf a}}$ becomes Fourier multipliers for
$\mathcal{M}_{p,q}(\R)$ for all $1<p<\infty$ and $1\leq q\leq \infty.$ More precisely, they proved
that
\begin{theorem}\cite{BGGO}\label{lp1} Let $\alpha>0$ and $\Omega=\{[\alpha n, \alpha(n+1)]:n\in \Z\}$.
For a bounded sequence ${\bf a}=\{a_n\}$ consider the function $\sigma_{\Omega,{\bf a}}$ as defined above.
Then $\sigma_{\Omega,{\bf a}}\in M(\mathcal{M}_{p,q}(\R))$ for all $1<p<\infty$ and $ 1\leq q\leq \infty.$
\end{theorem}
The authors have used Gabor frame characterization of modulation spaces in order to prove the above result.
But, we would like to remark that this theorem can be proved easily using the other definition
(definition~(\ref{def2})) of modulation space norm.
In~\cite{BGOR} B\'{e}nyi, Gr\"{o}chenig, Okoudjou, and Rogers generalized Theorem~\ref{lp1} to include other
collection of intervals and also pointed out the same remark. In particular, they proved that
\begin{theorem}\cite{BGOR}\label{lp2} Let $\{b_n\}_{n\in \Z}$ be an increasing sequence of real numbers such
that $\inf\limits_n |b_{n+1}-b_n|=\beta>0.$ Consider $\Omega=\{[b_n,b_{n+1}]:n\in \Z\}$. Then given a bounded
sequence ${\bf a}=\{a_n\}_{n\in \Z},$ the function $\sigma_{\Omega,{\bf a}}$  is a Fourier multiplier on modulation
spaces $\mathcal{M}_{p,q}(\R)$ for all $1<p<\infty$ and $ 1\leq q\leq \infty.$
\end{theorem}
Note that in Theorem~\ref{lp1} intervals are of equal lengths and moreover they are
translates of one single interval. Whereas, in Theorem~\ref{lp2}, authors
have a restriction namely lengths of intervals cannot be arbitrarily small. We observe that it is
not the lengths of intervals but the locations of intervals, which play a role to become Fourier multiplier on
modulation spaces. This is the underlying idea for many results presented in this section. In particular, we prove
the following generalization of Theorem~\ref{lp1}.
\begin{theorem} \label{re1} Let $\Omega=\{\omega_n~:n\in \Z\}$ be a collection of intervals such that for all
$n\in \Z,~\omega_n\subseteq [\alpha n, \alpha(n+1)]$ for some $\alpha>0.$ Then for all bounded sequences
${\bf a}=\{a_n\},$ the function $\sigma_{\Omega,{\bf a}}$ is a Fourier multiplier for $\mathcal{M}_{p,q}(\R)$
for all $1<p<\infty$ and $ 1\leq q\leq \infty.$
\end{theorem}
\noindent {\bf Proof:}
We first note that it suffices to prove Theorem~\ref{re1} with the assumption that  
$\omega_n\subseteq [n, n+1],~n\in \Z.$ This follows by using standard dilation arguments. 

Let $\phi$ be a Schwartz class function such that ${\it supp}~\phi \subseteq [-1,1]$ and $\sum_k \phi_k \equiv 1,$ where
$\phi_k(.)=\phi(.-k),$ i.e., $\phi_k$'s form a partition of unity. Note that ${\it supp}~\phi_k \subseteq [k-1,k+1].$
Hence, for fixed $k\in \Z,$ there exists at most two intervals $\omega_n,$ namely $\omega_{k-1}$ and $\omega_k,$
such that $\omega_n$ intersects with ${\it supp}~\phi_k.$ Thus we see that for all $k\in \Z,$  the function
$\phi_k \sigma \in M(L_p(\R))$ with norm independent of $k.$ This proves the desired result.
\qed

Next, we ask another natural question about Fourier multipliers on modulation spaces~: Does there exist
$\sigma\in L^{\infty}(\R^n)$ such that $\sigma\notin M(\mathcal{M}_{p,q}(\R^n)),~p\neq 2?$ We first answer this question
in general and then provide some concrete examples of bounded measurable functions which are not multipliers for
$\mathcal{M}_{p,q}(\R^n),~p\neq 2.$
\begin{proposition}\label{comp}
Let $\sigma\in L^{\infty}(\R^n)$ be such that $\sigma\notin M(L^p(\R^n)),~p\neq 2.$ Further, assume that $\sigma$
has compact support. Then, $\sigma\notin M(\mathcal{M}_{p,q}(\R^n)),~1\leq q\leq \infty.$
\end{proposition}
\noindent {\bf Proof:}
Let $p\neq 2$ and $\sigma\in L^{\infty}(\R^n)$ be a compactly supported function such that 
$\sigma\notin M(L^p(\R^n)).$ We are interested in proving that $\sigma\notin M(\mathcal{M}_{p,q}(\R^n)),~
1\leq q\leq \infty.$ Without loss of generality we may assume that 
${\it supp}~\sigma \subseteq[-\frac{1}{4},\frac{1}{4}]^n.$

Suppose on the contrary that $\sigma\in M(\mathcal{M}_{p,q}(\R^n)).$ We apply the definition of
Fourier multipliers on modulation spaces with a particular choice of partition of unity to arrive 
at a contradiction. 

Let $\phi\in \mathcal{S}(\R^n)$ be such that ${\it supp}~\phi \subseteq [-1,1],$ and $\sum_k \phi_k \equiv 1,$ where
$\phi_k(.)=\phi(.-k).$ In addition to this, we assume that  $\phi\equiv 1$ on $[-\frac{1}{4},\frac{1}{4}]^n.$ 
Since $\sigma\in M(\mathcal{M}_{p,q}(\R^n)),$ by using the definition of Fourier multipliers on modulation spaces, we
conclude that $\sigma \phi_k \in M(L^p(\R^n)),~k\in \Z^n.$ In particular, $\sigma \phi_0 \in M(L^p(\R^n)).$ But, 
our choice of the function $\phi,$ is such that we have $\sigma \phi_0=\sigma,$ which contradicts our hypothesis 
that $\sigma\notin M(L^p(\R^n)).$
This completes the proof. 
\qed
 
From the above proposition and the celebrated  ball multiplier result due to  C. Fefferman we have 
\begin{example}\label{ex2}
In~\cite{CF} C. Fefferman proved that characteristic function of the unit ball $B_1(0)$ in $\R^n, n\geq 2,$
is not an $L^p-$multiplier for $p\neq 2.$ We apply Proposition~\ref{comp} to conclude that $\chi_{B_1(0)}\notin M(\mathcal{M}_{p,q}(\R^n))$ for all $p\neq 2$ and $1\leq q\leq \infty.$
\end{example}
Example for $\R$ will follow from a  beautiful result  due to V. Lebedev and A. Olebski\^{i}~\cite{LO}. They proved that
\begin{theorem}\cite{LO}\label{lo}Let $E\subseteq \R^n$ be a measurable set such that $\chi_{_E}\in M(L^p(\R^n))$ for some $p\neq 2.$
Then $E$ is an open set upto a set of measure zero.
\end{theorem}
We apply this theorem together with Proposition~\ref{comp} to get the following example for modulation spaces.
\begin{example}Let $E\subseteq [0,1]$ be the Cantor set of positive measure. Then, 
$\chi_{_E} \notin M(\mathcal{M}_{p,q}(\R))$ for $p\neq 2.$
\end{example}

Now we will provide another example in this direction which will be useful later. We construct a partition of the interval $(0,1)$ into disjoint intervals such that it does not
give rise to an $L^p-$multiplier for $p\neq 2.$ Then as a consequence
of Proposition~\ref{comp} we can deduce that an arbitrary collection of disjoint intervals may not give rise to
Fourier multipliers on all modulation  spaces in the sense of Theorem~\ref{re1}.
\begin{example} \label{ex1}
Let $\omega=(0,1)$ and consider the dyadic partition of $\omega,$ i.e., $\omega=\cup_{n\geq 1} \omega_n,$
where $\omega_n=[2^{-n},2^{-n+1}).$  We further partition each interval $\omega_n$ into $2^n$ disjoint intervals of
equal lengths, i.e., for each $n\geq 1,$ we write $\omega_n =\cup_{m=1}^{2^n} \omega_{n,m},$ where
$\omega_{n,m}\subseteq \omega_n$ are disjoint intervals and $|\omega_{n,m}|=2^{-2n}$ for all $m=1,2,...,2^n.$

We claim that $\sigma(\xi)=\sum\limits_{n\geq 1}\sum\limits_{m=1}^{2^n}c_{n,m}\chi_{\omega_{n,m}}(\xi)$
is not an $L^p-$multiplier for $1<p\neq 2<\infty,$ where $\{c_{n,m}\}$ is an arbitrary bounded sequence of real
numbers. 

The proof is again by contradiction. Suppose that $\sigma \in M(L^p(\R)),~1<p\neq 2<\infty,$ for all
bounded sequences $\{c_{n,m}\}.$ We know that $\chi_{\omega_n}\in M(L^p(\R)),~1<p<\infty,$ with norm independent 
of $n.$ Further, we use the fact that product of two $L^p-$multipliers is again an $L^p-$multiplier and  
conclude that the functions $\sigma_n(\xi)=\sum\limits_{m=1}^{2^n}c_{n,m}\chi_{\omega_{n,m}}(\xi)\in M(L^p(\R)),$ for all 
$n\geq 1,$ with norm independent of $n.$
Note that for each fixed $n\geq 1,$ intervals $\omega_{n,m}$ are of equal lengths for all $m=1,2,...,2^n.$
Hence using dilation argument, we obtain that
functions of the form $\sum\limits_{m=1}^{2^n}c_{n,m}\chi_{[m,m+1)}(\xi)$ become Fourier multipliers for the
same $p$ with norm independent of $n.$ This is a contradiction to the fact that for an arbitrary bounded sequence
$\{c_m\},$ the function $\sum\limits_{m\in \Z}c_{m}\chi_{[m,m+1)}(\xi)$ is not an $L^p-$multiplier
unless $p=2.$

Since support of $\sigma$ is compact, as an application of Proposition~\ref{comp}, we conclude that 
$\sigma\notin M(\mathcal{M}_{p,q}(\R))$ for $p\neq 2.$ This
example also tells us that an arbitrary sequence of disjoint intervals may not give rise to Fourier multipliers
(in the sense of Theorem~\ref{re1}) on $\mathcal{M}_{p,q}(\R)$ for $p\neq 2.$
\end{example}

Next, we prove inclusion relations between multiplier spaces $M(\mathcal{M}_{p,q}(\R^n)).$ First observe that by
using duality and interpolation arguments, we have $M(\mathcal{M}_{p_1,q_1}(\R^n))\subseteq
M(\mathcal{M}_{p_2,q_2}(\R^n)),$ whenever $1\leq p_1\leq p_2\leq 2$ and $1\leq q_1,q_2\leq \infty.$ A similar
relation holds for dual exponents, i.e., $M(\mathcal{M}_{p_2,q_2}(\R^n))\subseteq M(\mathcal{M}_{p_1,q_1}(\R^n)),$
if $2\leq p_1\leq p_2\leq \infty$ and $1\leq q_1,q_2\leq \infty.$ Like $L^p-$multiplier spaces, the above mentioned inclusion
relations are strict and this is the content of next theorem.
\begin{theorem}\label{inclusion} Let $1\leq p_1<p_2\leq 2$ and
$1\leq q_1,q_2\leq \infty.$ Then, we have strict inclusion $M(\mathcal{M}_{p_1,q_1}(\R^n))\varsubsetneq
M(\mathcal{M}_{p_{2},q_{2}}(\R^n)).$
\end{theorem}
\noindent {\bf Proof:} In order to avoid certain notational inconvenience we only prove this theorem for $n=1.$
The higher dimensional analogue can be proved similarly. Let $p_1,p_2,q_1,q_2$ be as in theorem.

We need to find a function $\phi \in M(\mathcal{M}_{p_{2},q_{2}}(\R))$ such that
$\phi \notin M(\mathcal{M}_{p_1,q_1}(\R)).$ Note that if $p_1=1,$ then the classical Hilbert transform
provides the required function. More precisely, we know that for all $1<p<\infty, 1\leq q \leq \infty,$
${\it sgn}(\xi) \in M(\mathcal{M}_{p,q}(\R))$ and ${\it sgn}(\xi) \notin M(\mathcal{M}_{1,q}(\R)),$ as
$\mathcal{M}_1(\R))$ is invariant under Fourier transform. 

Assume that $1< p_1<p_2\leq 2.$ Observe that in order to prove the desired result, it is enough 
to show that there exists a compactly supported function 
$\sigma \in M(L^{p_2}(\R))$ such that $\sigma \notin M(L^{p_1}(\R)).$ This observation follows from Proposition~\ref{comp} and 
Property (4) of Proposition~\ref{prop}. The existance of the required function $\sigma,$ will be proved using  
some known classical transference results for $L^p-$multipliers due to de Leeuw~\cite{de} and Jodeit~\cite{Jod}.

Let $M(l_{p}(\Z))$ denote the space of
Fourier multipliers on $l_p(\Z).$ Notice that members of $M(l_{p}(\Z))$ are $1-$periodic functions. 
From the classical $L^p-$multiplier theory, we know that $M(l_{p_1}(\Z))\varsubsetneq M(l_{p_2}(\Z))$ whenever
$1\leq p_1<p_2\leq 2,$ i.e., there exists a function $m\in M(l_{p_2}(\Z))$ such that $m \notin M(l_{p_1}(\Z)).$
Let $m$ be such a function. Without loss of generality we may assume that
${\it supp~} m \subseteq [-\frac{1}{4},\frac{1}{4}].$
Let $m^{\sharp}$ denote the $1-$periodization of $m$ from the interval $[-\frac{1}{2},\frac{1}{2}).$
We apply de Leeuw's transference result (see~\cite{de})
about periodic multipliers to conclude that $m^{\sharp}\in M(L^{p_2}(\R)).$
Notice that $\chi_{[-\frac{1}{2},\frac{1}{2})}m=m.$ Hence the function $m^{\sharp}$ can also be thought of 
as $1-$periodization of the function
$\chi_{[-\frac{1}{2},\frac{1}{2})}m.$ Now, we use Jodeit's transference result~\cite{Jod} about periodization of
compactly supported multipliers to get that $\chi_{[-\frac{1}{2},\frac{1}{2})}m \in M(L^{p_2}(\R)).$ Since both
these transference results are if and only if type results, a repetition of previous arguments will lead 
us to the conclusion that $\chi_{[-\frac{1}{2},\frac{1}{2})}m \notin M(L^{p_1}(\R)).$
This completes the proof. 
\qed
\medskip

Similarly, we can prove that for $2\leq p_1<p_2\leq \infty$ and $1\leq q_1,q_2\leq \infty,$ we have
$M(\mathcal{M}_{p_2,q_2}(\R^n))\varsubsetneq M(\mathcal{M}_{p_1,q_1}(\R^n)).$

\section{\bf Littlewood-Paley Operators}\label{lp}
In this section we shall develop Littlewood-Paley theory for modulation spaces.
For this we shall need the notion of vector valued modulation spaces, which has a natural definition.
Since we are only interested in $l_2-$valued modulation spaces, we restrict ourselves to these spaces.
The $l_2-$valued modulation space, denoted by $\mathcal{M}_{p,q}(l_2),$ consists of sequences
$\{f_n\}$ of tempered distributions for which the following norm is finite:
\begin{eqnarray}\label{de1}
\left\|\{f_n\}\right\|_{\mathcal{M}_{p,q}(l_2)}= \left(\int_{\R^n} \left(\int_{\R^n}
\left(\sum_n|V_g(f_n)(x,\xi)|^2\right)^{\frac{p}{2}} dx\right)^{\frac{q}{p}}d\xi\right)^{\frac{1}{q}}.
\end{eqnarray}
The other definitions of modulation spaces can also be extended to $l_2-$ valued setting
in a similar fashion. For example, discrete version in vector valued setting takes the form:
\begin{eqnarray}\label{de2}
\left\|\{f_n\}\right\|_{\mathcal{M}_{p,q}(l_2)}= \left( \sum_k \left\|\left(\sum_n|T_{\phi_k}f_n|^2\right)^{\frac{1}{2}}\right\|^q_{L^p}\right)^{\frac{1}{q}}.
\end{eqnarray}
Given an intervals $\omega$ in $\R,$ let $S_{\omega}$ denote the multiplier operator associated with
symbol $\chi_{\omega}, $ i.e., $\widehat{S_{\omega} f}= \chi_{\omega} \hat f,~f\in \mathcal{S}(\R).$
Recall from previous section that $\chi_{\omega}\in M(\mathcal{M}_{p,q}(\R))$ for all $1<p<\infty$ and
$1\leq q\leq \infty.$
Moreover, operator norm $\|S_{\omega}\|$ is independent of the interval $\omega.$ Like
$L^p-$case, we consider Littlewood-Paley operators and investigate
their boundedness properties on modulation spaces. Let us first define Littlewood-Paley operators.
\begin{definition}[Littlewood-Paley Operator] \label{linearnonsmth:def} Let
$\Omega=\{\omega_n:n\in \Z\}$ be a collection of disjoint intervals in $\R$ and $S_{\omega_n}f$
be the multiplier operator defined as above. The Littlewood-Paley operator associated with
collection $\Omega$ is the $l_2-$valued operator $\Delta_{\Omega}:~f \rightarrow \{S_{\omega_n}f\},~f \in S(\R).$
\end{definition}
The theory of Littlewood-Paley operators on $L^p-$spaces is quite rich and it has many
beautiful applications in studying Fourier multiplier, characterizing important function spaces, etc.
Our concern in this paper is the boundedness properties of these Littlewood-Paley operators on modulation spaces.
For this we first recall some classical $L^p-$estimates for Littlewood-Paley operators.
\begin{theorem}\label{LP}
We have the following:
\begin{enumerate}
\item (Dyadic Littlewood-Paley~\cite{J}) Let $\Omega_d=\{\omega_n:n\in \Z\},$ where $\omega_n=(-2^{n+1},-2^n]\cup[2^n,2^{n+1}).$
Then for all $1<p<\infty,$ we have
\begin{eqnarray}\label{dyadic}\|\Delta_{\Omega_d}f\|_{L^p(\R)}\simeq \|f\|_{L^p(\R)},~\forall f\in \mathcal{S}(\R).
\end{eqnarray}
\item (Carleson~\cite{C}) If $\Omega=\{\omega_n=[n,n+1):n\in \Z\}.$ Then for all $2\leq p<\infty$
there exists a constant $C_p$ such that
\begin{eqnarray}\label{carleson}\|\Delta_{\Omega}f\|_{L^p(\R)}\leq C_p \|f\|_{L^p(\R)},~\forall f\in \mathcal{S}(\R).
\end{eqnarray}
\item (Rubio de Francia~\cite{R1}) Let $\Omega=\{\omega_n:n\in \Z\}$ be an arbitrary collection of
disjoint intervals. Then for all $2\leq p<\infty$ there exists a constant $C_p$ such that
\begin{eqnarray}\label{rubio}\|\Delta_{\Omega}f\|_{L^p(\R)}\leq C_p \|f\|_{L^p(\R)},~\forall f\in \mathcal{S}(\R).
\end{eqnarray}
\end{enumerate}
Moreover, $p\geq 2$ is a necessary condition in estimates~(\ref{carleson}) and~(\ref{rubio}).
\end{theorem}
In this section we study an analogue of Theorem~\ref{LP} in the context of modulation spaces. We would like
to remark here that we need to produce different (from classical $L^p-$case) arguments to prove analogue of
Theorem~\ref{LP} for modulation spaces. We also would like to mention that at many places
we (without mentioning it) will be dealing with only finite sequence of functions and operators so that all the steps are justified. Since, the estimates we obtain do not depend on sizes of sequences under consideration, we get the desired result using standard
limiting arguments. We first prove the following vector valued inequality~:
\begin{theorem} \label{cor2}Let $\{\omega_n\}_{n\in \Z}$ be a sequence of intervals in $\R$ and $S_{\omega_n}$
be the multiplier operator associated with symbol $\chi_{\omega_n}$. Then for $1<p<\infty$ and $1\leq q \leq \infty,$
there exists a constant $C_{p,q}$ such that for all sequences $\{f_n\}$, we have
\begin{eqnarray}\label{l2}
\left\|\{S_{\omega_n}f_n\}\right\|_{\mathcal{M}_{p,q}(l_2)}
\leq C_{p,q} \left\|\{f_n\}\right\|_{\mathcal{M}_{p,q}(l_2)}.
\end{eqnarray}
\end{theorem}

\noindent {\bf Proof:}
In order to prove this result we use its classical $L^p$ analogue, which says that for a sequence of intervals
$\omega_n$ and $1<p<\infty,$ there exists a constant $C_p$ such that

\begin{eqnarray}\label{lpextn2}
\left\|\{S_{\omega_n}f_n\}\right\|_{L^p(l_2)}
\leq C_p \left\|\{f_n\}\right\|_{L^p(l_2)},
\end{eqnarray}
where $\|\{f_n\}\|_{L^p(l_2)}=\|(\sum_n|f_n|^2)^{\frac{1}{2}}\|_{L^p(\R)}.$
Consider,
\begin{eqnarray*}
\left\|\{S_{\omega_n}f_n\}\right\|_{\mathcal{M}_{p,q}(l_2)}&=&
\left(\sum_k\left\|\left(\sum_n|T_{\phi_k}S_{\omega_n}f_n|^2\right)^{\frac{1}{2}}\right\|^q_{L^p}\right)^{\frac{1}{q}}\\
&=&\left(\sum_k\left\|\left(\sum_n|S_{\omega_n}T_{\phi_k}f_n|^2\right)^{\frac{1}{2}}\right\|^q_{L^p}\right)^{\frac{1}{q}}\\
&\leq&C_p\left(\sum_k\left\|\left(\sum_n|T_{\phi_k}f_n|^2\right)^{\frac{1}{2}}\right\|^q_{L^p}\right)^{\frac{1}{q}}\\
&=& C_p \left\|\{f_n\}\right\|_{\mathcal{M}_{p,q}(l_2)}.
\end{eqnarray*}
Here we have used vector valued inequality~(\ref{lpextn2}) together with the fact that operator $S_{\omega_n}$
and $T_{\phi_k}$ commute. This completes the proof.\qed
\begin{remark} We would like to remark that inequality (\ref{l2}) holds true even if we replace
$l_2$ with $l_r$ for any $1<r<\infty$ as its classical variant (\ref{lpextn2}) is known to be true for $l_r$
for all $1<r<\infty$ and $1<p<\infty.$
\end{remark}

\begin{theorem}\label{re3}
Let $\Omega=\{\omega_n:n\in \Z\}$ be a collection of intervals in $\R.$ If
$\sigma_{\Omega,{\bf a}}=\sum_n a_n \chi_{\omega_n}\in M(\mathcal{M}_{p,q}(\R))$ for all bounded sequences ${\bf a}=\{a_n\},$
with norm bounded by a constant multiple of $\|{\bf a}\|_{l_{\infty}}.$
Then, the Littlewood-Paley operator $\Delta_{\Omega}$ associated with the collection $\Omega$
is bounded on $\mathcal{M}_{p,q}(\R),$ i.e., we have the following:
\begin{eqnarray}
\left\|\{S_{\omega_n}f\}\right\|_{\mathcal{M}_{p,q}(l_2)}
\leq C_{p,q} \|f\|_{\mathcal{M}_{p,q}(\R)}.
\end{eqnarray}
\end{theorem}
\noindent {\bf Proof:}
Let $r_n$ denote the sequence of Radamacher functions. As an application of Khintchine's inequality
we know that for $0<p<\infty$ and sequence of complex numbers $b_n,$ we have
\begin{eqnarray*}\left(\sum\limits_n |b_n|^2\right)^{\frac{1}{2}} \simeq
\left(\int_0^1 |\sum\limits_n b_n r_n(t)|^p dt\right)^{\frac{1}{p}}.
\end{eqnarray*}
With the help of above inequality we linearize Littlewood-Paley operators, which is a standard
technique to deal with such objects. But, we would like to point out that the proof of this theorem is not
as straight forward as in the classical $L^p-$case. We need to consider two cases $p\leq q$ and
$p\geq q$ separately.

{\bf Case 1: $p\geq q$.} Consider
\begin{eqnarray*}
\left\|\Delta_{\Omega}f\right\|^q_{\mathcal{M}_{p,q}(l_2)}&=&
\left\|\{S_n f\}\right\|^q_{\mathcal{M}_{p,q}(l_2)}\\
&=&\int_\R \left(\int_\R \left(\sum_n|V_g(S_nf)(x,\xi)|^2\right)^{\frac{p}{2}} dx\right)^{\frac{q}{p}}d\xi\\
&& (\text {use Khintchine's inequality for exponent}~ q)\\
&\leq& C_q \int_\R \left(\int_\R \left(\int_0^1|\sum_n r_n(t)V_g(S_nf)(x,\xi)|^q dt\right)^{\frac{p}{q}} dx\right)^{\frac{q}{p}}d\xi\\
&=& C_q \int_\R \left(\int_\R \left(\int_0^1|V_g (\sum_n r_n(t)S_n f)(x,\xi)|^q dt\right)^{\frac{p}{q}} dx\right)^{\frac{q}{p}}d\xi\\
&& (\text {apply Minkowski's integral inequality with exponent~} p/q)\\
&\leq& C_q \int_\R \int_0^1 \left(\int_\R |V_g (H_{\Omega, \{r_n(t)\}} f)(x,\xi)|^p dx\right)^{\frac{q}{p}} dt d\xi\\
&=& C_q \int_0^1 \|H_{\Omega, \{r_n(t)\}}f\|_{\mathcal{M}_{p,q}}^q dt\\
&=& C'_q \|f\|^q_{\mathcal{M}_{p,q}(\R)}.
\end{eqnarray*}
where $H_{\Omega, \{r_n(t)\}}$ is the Fourier multiplier operator associated with the symbol
$\sum_n r_n(t)\chi_{\omega_n}$ and we have used that the operator norm of $H_{\Omega, \{r_n(t)\}}$
is uniformly bounded in $t.$

{\bf Case 2: $p\leq q$.} This time we use Khintchine's inequality for exponent $p$ as follows~:
\begin{eqnarray*}
\left\|\Delta_{\Omega}f\right\|^q_{\mathcal{M}_{p,q}(l_2)}
&=&\int_\R \left(\int_\R \left(\sum_n|V_g(S_nf)(x,\xi)|^2\right)^{\frac{p}{2}} dx\right)^{\frac{q}{p}}d\xi \\
&& (\text {use Khintchine's inequality for exponent}~ p)\\
&\leq& C_p \int_\R \left(\int_\R \int_0^1|\sum_n r_n(t)V_g(S_nf)(x,\xi)|^p dt dx\right)^{\frac{q}{p}}d\xi\\
&=& C_p \int_\R \left(\int_0^1\int_\R|V_g (\sum_n r_n(t)S_n f)(x,\xi)|^p dx dt\right)^{\frac{q}{p}}d\xi\\
&\leq& C_p \int_\R \int_0^1 \left(\int_\R |V_g (H_{\Omega, \{r_n(t)\}} f)(x,\xi)|^p dx\right)^{\frac{q}{p}} dt d\xi\\
&=& C_p \int_0^1 \|H_{\Omega, \{r_n(t)\}}f\|_{\mathcal{M}_{p,q}}^q dt\\
&\leq& C'_p \|f\|^q_{\mathcal{M}_{p,q}}.
\end{eqnarray*}
This completes the proof.
\qed

As an immediate application of the above theorem together with Theorems~\ref{lp2} and~\ref{re1},
we conclude the following result for Littlewood-Paley operators on modulation spaces:
\begin{theorem}
Let $\Omega=\{\omega_n:n\in \Z\}$ be a collection of disjoint intervals. Assume that intervals are
either dyadic or satisfy hypothesis of any of Theorems~\ref{lp2} and~\ref{re1}. Then the associated
Littlewood-Paley operator $\Delta_{\Omega}$ maps  $\mathcal{M}_{p,q}(\R)$ into $\mathcal{M}_{p,q}(l_2)$
for all $1<p<\infty$ and $1\leq q \leq \infty.$
\end{theorem}
Note that in the above theorem we have obtained modulation spaces analogues of dyadic and Carleson's Littelwood-Paley 
results given by Theorem~\ref{LP}(1) and (2) respectively. Now we proceed to prove the analogue of Rubio de-Francia's 
Littlewood-Paley result in the context of modulation spaces. Recall Example~\ref{ex1}, where we have proved that an arbitrary collections of disjoint intervals may not
give rise to Fourier multipliers on modulation spaces $\mathcal{M}_{p,q}(\R)$ for $p\neq 2.$
But, we will see that for an arbitrary collection of disjoint intervals an analogue of Theorem~\ref{LP}(3)
for modulation spaces holds. More precisely, we have
\begin{theorem}\label{re4}
Let $\Omega=\{\omega_n:n\in \Z\}$ be a collection of disjoint intervals. Then the associated
Littlewood-Paley operator $\Delta_{\Omega}$ is bounded from  $\mathcal{M}_{p,q}(\R)$
into $\mathcal{M}_{p,q}(l_2)$ for all $2\leq p<\infty$ and $1\leq q \leq \infty.$
\end{theorem}
\noindent {\bf Proof:}
We shall use the discrete version (see~(\ref{de2})) of definition for modulation space norm. For
all $f\in \mathcal{S}(\R)$ and $2\leq p<\infty$ and $1\leq q \leq \infty,$ we need to prove that
\begin{eqnarray*}
\left\|\Delta_{\Omega}f\right\|_{\mathcal{M}_{p,q}(l_2)}&=&
\left\|\{S_n f\}\right\|_{\mathcal{M}_{p,q}(l_2)}\\
&=& \left( \sum_k \left\|\left(\sum_n|T_{\phi_k}S_n f|^2\right)^{\frac{1}{2}}\right\|^q_{L^p}\right)^{\frac{1}{q}}\\
&\leq& C_{p,q}\|f\|_{\mathcal{M}_{p,q}(l_2)}
\end{eqnarray*}
Observe that for all $n,k\in \Z,$ we have $T_{\phi_k}S_n f=S_n T_{\phi_k}f.$ Hence using Rubio de
Francia's Littlewood-Paley Theorem~\ref{LP}(3), we have
\begin{eqnarray*}
\left\|\left(\sum_n|T_{\phi_k}S_n f|^2\right)^{\frac{1}{2}}\right\|_{L^p}&=& \left\|\left(\sum_n|S_nT_{\phi_k} f|^2
\right)^{\frac{1}{2}}\right\|_{L^p}\\
&\leq& C_p \left\|T_{\phi_k} f\right\|_{L^p},~2\leq p<\infty.
\end{eqnarray*}
Substituting this estimate in above we get that
\begin{eqnarray*}
\left\|\Delta_{\Omega}f\right\|_{\mathcal{M}_{p,q}(l_2)}
&\leq& \left( \sum_k C_p\left\|T_{\phi_k} f\right\|^q_{L^p}\right)^{\frac{1}{q}}\\
&=& C_p\|f\|_{\mathcal{M}_{p,q}(\R)}.
\end{eqnarray*}
This proves boundedness of the Littlewood-Paley operator $\Delta_{\Omega}$ and hence the proof is complete.
\qed

Next, we prove that $p\geq 2$ is a necessary condition in the above theorem. We would like to remark here
that in the case of classical $L^p-$spaces, the necessity of $p\geq 2$ in Theorem~\ref{LP}(2)
(and hence in Theorem~\ref{LP}(3))
is proved by getting an estimate for the square function associated with the sequence $\{[n,n+1)\}$
for a particular choice of  function. But, this does not work in the case of modulation spaces as we have
seen in Theorem~\ref{lp1} that this sequence of intervals even gives rise to Fourier multipliers on all
modulation spaces. But interestingly, $p\geq 2$ still remains a necessary condition
in Theorem~\ref{re4}, unlike Carleson's analogue. In order to prove this, we require the following proposition.
\begin{proposition}\label{necepro} Let $\Omega=\{\omega_n:n\in \Z\}$ be a collection of disjoint intervals. Let
$1<p\leq 2$ and $1\leq q\leq \infty.$ Assume that the associated Littlewood-Paley operator $\Delta_{\Omega}$ is bounded
from $\mathcal{M}_{p,q}(\R)$ into $\mathcal{M}_{p,q}(l_2).$
Then, $\sigma_{\Omega,{\bf a}}(\xi)=\sum_n a_n \chi_{\omega_n}(\xi)\in M(\mathcal{M}_{p,q}(\R)),$ where
${\bf a}=\{a_n\}\in l_{\infty}.$ Moreover, the norm is bounded by a constant multiple of
$\|{\bf a}\|_{l_{\infty}}\|\Delta_{\Omega}\|.$
\end{proposition}
\noindent {\bf Proof:}
Let $\Omega=\{\omega_n:n\in \Z\}$ be a given collection of disjoint intervals in $\R$ and
${\bf a}=\{a_n\}\in l_{\infty}$ be a given bounded sequence. Let $H_{\Omega,{\bf a}}$ denote the
multiplier operator associated with the symbol $\sigma_{\Omega,{\bf a}}=\sum_n a_n \chi_{\omega_n}.$
We need to prove that for all $f\in \mathcal{S}(\R),$
\begin{eqnarray*}
\|H_{\Omega,{\bf a}}f \|_{\mathcal{M}_{p,q}}
&=& \left( \sum_k \|T_{\phi_k}H_{\Omega,{\bf a}}f\|^q_{L^p}\right)^{\frac{1}{q}}\\
&\leq& C_{p,q}\|f\|_{\mathcal{M}_{p,q}}
\end{eqnarray*}
For $k\in \Z,$ consider,
\begin{eqnarray*}
|\langle T_{\phi_k}H_{\Omega,{\bf a}}f, g\rangle|
&=& |\langle H_{\Omega,{\bf a}}T_{\phi_k}f, g\rangle|\\
&=& |\int_\R \sum\limits_n a_n S_{\omega_n}T_{\phi_k}f(x)S_{\omega_n}g(x)|dx\\
&\leq & \|{\bf a}\|_{l_\infty} \int_R \left(\sum_n|S_{\omega_n} T_{\phi_k}f(x)|^2\right)^{\frac{1}{2}}
\left(\sum_n|S_{\omega_n} g(x)|^2\right)^{\frac{1}{2}}dx\\
&\leq & \|{\bf a}\|_{l_\infty} \left\|\left(\sum_n|S_{\omega_n} T_{\phi_k}f|^2\right)^{\frac{1}{2}} \right\|_{L^p}
\left\|\left(\sum_n|S_{\omega_n} g|^2\right)^{\frac{1}{2}}\right\|_{L^{p'}}\\
& & (\text {As}~ p'\geq 2, \text {use Rubio de Francia Littlewood-Paley Theorem}~\ref{LP}(3))\\
&\leq & C_p \|{\bf a}\|_{l_\infty} \left\|\left(\sum_n|S_{\omega_n} T_{\phi_k}f|^2\right)^{\frac{1}{2}}
\right\|_{L^p} \|g\|_{L^{p'}}.
\end{eqnarray*}
Since this holds for all $g\in L^{p'}(\R),$ we get
$$\|T_{\phi_k}H_{\Omega,{\bf a}}f\|_{L^{p}}\leq C_p \|{\bf a}\|_{l_\infty} \left\|\left(\sum_n|S_{\omega_n} T_{\phi_k}f|^2\right)^{\frac{1}{2}}
\right\|_{L^p}.$$ Hence we have,
\begin{eqnarray*}
\|H_{\Omega,{\bf a}}f \|_{\mathcal{M}_{p,q}}
&=& C_p \|{\bf a}\|_{l_\infty} \left( \sum_k \left\|\left(\sum_n|S_{\omega_n} T_{\phi_k}f|^2\right)^{\frac{1}{2}}
\right\|_{L^p}^q\right)^{\frac{1}{q}}\\
&=& C_p \|{\bf a}\|_{l_\infty} \|\Delta_{\Omega}f\|_{\mathcal{M}_{p,q}}\\
&\leq& C_{p,q} \|{\bf a}\|_{l_\infty}\|f\|_{\mathcal{M}_{p,q}}.
\end{eqnarray*}
Here we have used the assumption that $\Delta_{\Omega}$ maps $\mathcal{M}_{p,q}(\R)$ into $\mathcal{M}_{p,q}(l_2)$
and
this finishes the proof.
\qed

\begin{corollary}\label{nece}
$p\geq 2$ is a necessary condition in Theorem~\ref{re4}.
\end{corollary}
\noindent {\bf Proof:}
The proof follows by considering the collection of intervals discussed in Example~\ref{ex1} with Proposition~\ref{necepro}.
Let $1<p<2$ and $1\leq q\leq \infty.$  Now suppose on the contrary that for this range of $p$ and $q,$
Theorem~\ref{re4} holds for all collections of disjoint intervals. Hence, in particular, it holds for the collection of
intervals described in Example~\ref{ex1}. Let us denote that collection of intervals as
$\Omega=\{\omega_n\}.$ With our assumption we get that the associated Littlewood-paley operator $\Delta_{\Omega}$
is bounded from $\mathcal{M}_{p,q}(\R)$ into $\mathcal{M}_{p,q}(l_2).$
As a consequence of Proposition~\ref{necepro}, we see that
for all bounded sequences ${\bf a}=\{a_n\},$ the function $\sigma_{\Omega,{\bf a}}=\sum_n a_n \chi_{\omega_n} \in
M(\mathcal{M}_{p,q}(\R)).$ But, this contradicts the fact (see Example~\ref{ex1}) that $\sigma_{\Omega,{\bf a}}$
may not be a Fourier multiplier for $\mathcal{M}_{p,q}(\R),$ if $p\neq 2.$
Thus we arrive at a contradiction.
\qed

\section{\bf A theorem of Marcinkiewicz and Zygmund for modulation spaces}\label{mz}
It is a well known classical result due to Marcinkiewicz and Zygmund that any bounded linear operator from
$L^p(\R^n)$ into itself admits an $l_2-$valued extension. In this section we shall prove
an analogue of this result for bounded linear operators on modulation spaces. We would like to remark
that proof of this result is quite different from its classical variant as there are two parameters
$p$ and $q$ in case of modulation spaces.
\begin{theorem} \label{l2ext} Let $1\leq p,q\leq \infty.$  Assume that $T$ is a bounded linear
operator from $\mathcal{M}_{p,q}(\R^n)$ into itself. Then, $T$ admits an $l_2-$valued
bounded extension. Moreover, the operator norm is bounded by a constant multiple of $\|T\|,$ where $\|T\|$ is
the operator norm of $T$ on $\mathcal{M}_{p,q}(\R^n).$
\end{theorem}
\noindent {\bf Proof:}
Let $\{f_n\}\in \mathcal{M}_{p,q}(l_2).$ We may assume that it is a finite sequence. We need to prove the following $l_2-$valued
estimate for the operator $T$~:
\begin{eqnarray}\label{extn} \left\|\{Tf_n\}\right\|_{\mathcal{M}_{p,q}(l_2)}&\leq& C_{p,q}\|T\| \left\|\{f_n\}
\right\|_{\mathcal{M}_{p,q}(l_2)}.
\end{eqnarray}
We need to consider the cases $p\geq q$ and $p\leq q$ separately.

{\bf Case 1.~$p\geq q~:$}
\medskip
We linearize the $l_2-$norm of a sequence with the help of Radamacher functions $r_n$ using the
Khintchine's inequality. Consider,
\begin{eqnarray*}
\left\|\{Tf_n\}\right\|^q_{\mathcal{M}_{p,q}(l_2)}&=& \int_\R \left(\int_\R \left(\sum_n|V_g(Tf_n)(x,\xi)|^2\right)^{\frac{p}{2}} dx\right)^{\frac{q}{p}}d\xi\\
&& (\text {use Khintchine's inequality with exponent}~ q)\\
&\leq& C_q \int_\R \left(\int_\R \left(\int_0^1|\sum_n r_n(t)V_g(Tf_n)(x,\xi)|^q dt\right)^{\frac{p}{q}} dx\right)^{\frac{q}{p}}d\xi\\
&=& C_q \int_\R \left(\int_\R \left(\int_0^1|V_g (T(\sum_n r_n(t)f_n))(x,\xi)|^q dt\right)^{\frac{p}{q}} dx\right)^{\frac{q}{p}}d\xi\\
&& (\text {as}~ p\geq q,~ \text {apply Minkowski's integral inequality with}~ p/q)\\
&\leq& C_q \int_\R \int_0^1 \left(\int_\R |V_g (T(\sum_n r_n(t)f_n))(x,\xi)|^p dx\right)^{\frac{q}{p}} dt d\xi\\
&=& C_q \int_0^1 \|T(\sum_n r_n(t)f_n)\|_{\mathcal{M}_{p,q}}^q dt\\
&\leq& C_q \|T\|^q \int_0^1 \|\sum_n r_n(t)f_n\|_{\mathcal{M}_{p,q}}^q dt\\
&=& C_q \|T\|^q\int_\R \int_0^1\left(\int_\R |\sum_n r_n(t)V_gf_n(x,\xi)|^p dx\right)^{\frac{q}{p}} dt d\xi\\
&\leq& C_q \|T\|^q\int_\R \left(\int_0^1\int_\R |\sum_n r_n(t)V_gf_n(x,\xi)|^p dx dt\right)^{\frac{q}{p}} d\xi\\
&& (\text {now use Khintchine's inequality with exponent}~ p)\\
&\leq& C_{p,q} \|T\|^q \left\|\{f_n\}\right\|^q_{\mathcal{M}_{p,q}(l_2)}.
\end{eqnarray*}
\medskip

{\bf Case 2.~$p\leq q~:$} By using Khintchine's inequality with exponent $p,$ we have
\begin{eqnarray*}
\left\|\{Tf_n\}\right\|^q_{\mathcal{M}_{p,q}(l_2)}
&\leq& C_{p} \int_\R \left(\int_\R \int_0^1|V_g (T(\sum_n r_n(t)f_n))(x,\xi)|^{p} dt dx\right)^{\frac{q}{p}}d\xi\\
&\leq& C_{p} \int_\R \int_0^1 \left(\int_\R |V_g (T(\sum_n r_n(t)f_n))(x,\xi)|^p dx\right)^{\frac{q}{p}} dt d\xi\\
&=& C_{p} \int_0^1 \|T(\sum_n r_n(t)f_n)\|_{\mathcal{M}_{p,q}}^q dt\\
&\leq& C_{p} \|T\|^q \int_0^1 \|\sum_n r_n(t)f_n\|_{\mathcal{M}_{p,q}}^q dt\\
&=& C_{p} \|T\|^q \int_\R \left(\int_0^1\left(\int_\R |\sum_n r_n(t)V_gf_n(x,\xi)|^p dx\right)^{\frac{q}{p}} dt\right)^{\frac{p}{q}.\frac{q}{p}} d\xi\\
&& (\text {apply Minkowski's integral inequality with exponent}~q/p)\\
&\leq& C_{p} \|T\|^q\int_\R \left(\int_\R \left(\int_0^1|\sum_n r_n(t)V_gf_n(x,\xi)|^q  dt\right)^{\frac{p}{q}} dx\right)^{\frac{q}{p}}d\xi\\
&& (\text {now use Khintchine's inequality with exponent}~ q)\\
&\leq& C_{p,q} \|T\|^q \left\|\{f_n\}\right\|^q_{\mathcal{M}_{p,q}(l_2)}.
\end{eqnarray*}
This completes the proof of Theorem~\ref{l2ext}.
\qed

%
%

\end{document}